\documentclass[11pt]{amsart}
\usepackage{amssymb,latexsym,epsfig,color}

\newcommand{\Z}{{\mathbb Z}}
\newcommand{\N}{{\mathbb N}}
\newcommand{\R}{{\mathbb R}}

\newcommand{\A}{{\mathcal A}}

\newcommand{\Set}{{\mathcal S}}

\newcommand{\G}{{\mathcal G}}

\newtheorem{theorem}{Theorem}[section]
\newtheorem{proposition}[theorem]{Proposition}
\newtheorem{corollary}[theorem]{Corollary}
\newtheorem{lemma}[theorem]{Lemma}

\newtheorem{problem}[theorem]{Problem}

\theoremstyle{definition}

\newtheorem{remark}[theorem]{Remark}
\newtheorem{example}[theorem]{Example}

\title [Alternatives for TDI] {Alternatives for Testing  Total Dual
Integrality}
\author{Edwin O'Shea and Andr{\' a}s Seb{\H o} }
\thanks{
E.O'S. was supported by a Fulbright grant and by
NSF grants DMS-9983797 and DMS-0401047.
A.S. was supported by the
``Marie Curie Training Network'' ADONET of the European Community.}
\address{Departamento de Matem\'aticas,
Centro de Investigaci\'on y de Estudios Avanzados del IPN,
Apartado Postal 14--740, 07000 M\'exico D.F., M{\'E}XICO and
CNRS, Laboratoire G-SCOP,
 46, Avenue F\'elix Viallet, 38000 Grenoble 38031 Grenoble, Cedex 1, FRANCE}
\email{edwin@math.cinvestav.mx, Andras.Sebo@g-scop.inpg.fr}
\date{\today}

\begin{document}

\maketitle

\begin{abstract}
In this paper we provide characterizing properties of TDI systems,
among others the following: a system of linear inequalities is TDI
if and only if its coefficient vectors form a Hilbert basis, and
there exists a test-set for the system's dual integer programs where
all test vectors have positive entries equal to $1$. Reformulations
of this provide relations between computational algebra and integer
programming and they contain Applegate, Cook and McCormick's
sufficient condition for the TDI property and  Sturmfels' theorem
relating toric initial ideals generated by square-free monomials to
unimodular triangulations. We also study the theoretical and
practical efficiency and limits of the characterizations of the TDI
property presented here.

In the particular case of set packing polyhedra our results
correspond to endowing the weak perfect graph theorem with an
additional, computationally interesting, geometric feature: the
normal fan of the stable set polytope of a perfect graph can be
refined into a regular triangulation consisting only of unimodular
cones.

\end{abstract}

\section{Introduction}

A restricted draft concerning an earlier stage of this research has been
reported in the IPCO 2007 conference proceedings \cite{os}.

Let $A = [{\bf a}_1 \, {\bf a}_2 \cdots {\bf a}_n ] \in \Z^{d \times
n}$ and assume that $A$ has rank $d$. With an abuse of notation the
ordered set of vectors consisting of the columns of $A$ will also be
denoted by $A$. For every $\sigma \subseteq [n] := \{ 1, \ldots, n
\}$ we have the $d \times |\sigma|$ matrix $A_\sigma$ given by the
columns of $A$ indexed by $\sigma$. Let $\textup{cone}(A)$, $\Z A$
and $\N A$ denote the non-negative real, integer and non-negative
integer span of $A$ respectively and assume that $\Z A = \Z^d$.

Fixing ${\bf c} \in {\R}^n$, for each ${\bf b} \in \R^d$ the
{\em linear program} (or {\em primal program})
$\textup{LP}_{A,{\bf c}}({\bf b})$ and its
{\em dual program} $\textup{DP}_{A,{\bf c}}({\bf b})$ are defined by
$$
\textup{LP}_{A,{\bf c}}({\bf b}) \, := \,
                    \textup{minimize} \, \{ \, {\bf c} \cdot {\bf x} \, :
            \, A{\bf x} \, = \, {\bf b}, \,
               {\bf x} \geq {\bf 0} \, \}
$$
and
$
\textup{DP}_{A,{\bf c}}({\bf b}) \, := \,
                    \textup{maximize} \, \{ \, {\bf y} \cdot {\bf b} \, :
            \, {\bf y}A \,  \leq \, {\bf c} \, \}.
$ Let $P_{\bf b}$ and $Q_{\bf c}$ denote the feasible regions of
$\textup{LP}_{A,{\bf c}}({\bf b})$ and $\textup{DP}_{A,{\bf c}}({\bf
b})$ respectively. Note that the linear program $\textup{LP}_{A,{\bf
c}}({\bf b})$ is feasible if and only if ${\bf b} \in
\textup{cone}(A)$. We refer to Schrijver \cite{schrijver} for basic
terminology, facts and notations about linear programming.

The corresponding {\em integer program} is defined as
$$
\textup{IP}_{A,{\bf c}}({\bf b}) \, := \,
                    \textup{minimize} \, \{ \, {\bf c} \cdot {\bf x} \, :
            \, A{\bf x} \, = \, {\bf b}, \,
               {\bf x} \in \N^n \, \}.
$$
We say that ${\bf c} \in \R^n$ is {\em generic} for $A$ if the integer
program $\textup{IP}_{A, {\bf c}}({\bf b})$ has a unique optimal solution
for all ${\bf b} \in \N A$. If ${\bf c}$ is generic
then $Q_{\bf c} \ne \emptyset$ and each linear program
$\textup{LP}_{A, {\bf c}}({\bf b})$ also has a unique optimal
solution for all ${\bf b} \in \textup{cone}(A)$ but the converse is
not true in general. However, the converse clearly holds for TDI
systems: the system ${\bf y}A \leq {\bf c}$ is {\em totally dual
integral} {\em (TDI} ) if $\textup{LP}_{A,{\bf c}}({\bf b})$ has an
integer optimal solution ${\bf x} \in \N^n$ for each ${\bf b} \in
\textup{cone}(A) \cap \Z^d $. In other words, the system ${\bf y}A
\leq {\bf c}$ defining $Q_{\bf c}\ne\emptyset$ is TDI exactly if the
optimal values of $\textup{LP}_{A,{\bf c}}({\bf b})$ and of
$\textup{IP}_{A,{\bf c}}({\bf b})$ coincide for all ${\bf b} \in
\textup{cone}(A) \cap \Z^d$. (If $Q_{\bf c} \ne \emptyset$, LP and
DP are both feasible and bounded.) This is a slight twist of
notation when compared to habits in combinatorial optimization: we
defined the TDI property for the dual problem, in order to be in
accordance with notations in computational algebra.

Totally dual integral (TDI) systems of linear inequalities play a
central role in combinatorial optimization. The recognition of TDI
systems has been recently proved to be coNP-complete by Ding, Feng
and Zang \cite{chinese} and this result has been sharpened to the
recognition of explicitly given systems with only $0-1$ coefficient
vectors, and where the defined polyhedron has exactly one vertex
(Hilbert basis testing by Pap \cite{juli}). Graph
theory results of Chudnovsky, Cornu{\' e}jols, Liu, Seymour and
Vu\v{s}kovi\'c \cite{recognition} allow one to recognize TDI systems
with $0-1$ coefficient matrices  and right hand sides. However,
solving the corresponding dual pair of integer linear programs
(including the coloration of perfect graphs) in polynomial time with
combinatorial algorithms remains open even in this special case. The
fixed dimension case has been solved long ago \cite{cls}, whereas the
fixed codimension case only recently \cite{ess}, \cite{hs}, the first
using a generalization of integer programming, the second computer
algebra,   both starting from a characterization of Hilbert bases in
\cite{sebo}.

A particular case where the recognition of  TDI
systems is still open, is the case of generic systems
(Problem~\ref{prob:unimod}), which is slightly more general than the
perfectness test (detecting perfection in a graph), and could be a
possible start for an alternative, simpler algorithm for the latter.
This paper wishes to contribute to testing TDI 
in cases that occur in integer programming and combinatorial
optimization, and which usually do not belong to the extremities 
that have been understood so far.

In Section~\ref{sec:TDI}, characterizing properties of TDI systems
are provided. Some of these properties involve tools from
combinatorial optimization, some others from computational algebra.
Section~\ref{sec:perfect} specializes these results to integral set
packing polytopes. Finally, Section~\ref{sec:computation} exhibits
the possible alternatives and their relative efficiency for
recognizing TDI systems. The remainder of this introduction is devoted
to explaining the main results and providing some of the necessary background.

A collection of subsets $\{\sigma_1, \ldots, \sigma_t\}$ of
$[n]$ will be called a {\em regular subdivision} of $A$ if there
exists ${\bf c} \in \R^n$, and ${\bf z}_1, \ldots, {\bf z}_t \in
\R^d$, such that ${\bf z}_i \cdot {\bf a}_j = c_j$ for all
 $j \in \sigma_i$ and ${\bf z}_i \cdot {\bf a}_j < c_j$
for all $j \notin \sigma_i$. The sets $\sigma_1, \ldots, \sigma_t$
are called the {\em cells} of the regular subdivision and the
regular subdivision is denoted by $\Delta_{\bf c}(A) = \{ \sigma_1,
\ldots, \sigma_t \}$ or simply $\Delta_{\bf c}$ when $A$ is unambiguous.

Equivalently, regular subdivisions are simply capturing {\em
complementary slackness} from linear programming. Namely, a feasible
solution to $\textup{LP}_{A,{\bf c}}({\bf b})$ is optimal if and
only if the support of the feasible solution is a subset of some
cell of $\Delta_{\bf c}$. Geometrically, $\Delta_{\bf c}$ can be
thought of as a partition of $\textup{cone}(A)$ by the inclusionwise
maximal ones among the cones $\textup{cone}(A_{\sigma_1}), \ldots,
\textup{cone}(A_{\sigma_t})$; each such cone is generated by the
normal vectors of defining inequalities of faces of $Q_{\bf c}$,
each maximal cell indexes the set of normal vectors of inequalities
satisfied with equality by a vertex (or minimal face) of $Q_{\bf
c}$. The ${\bf z}_1, \ldots, {\bf z}_t$ above are these
vertices and so the regular subdivision $\Delta_{\bf c}$
is geometrically realized as the {\em normal fan} of $Q_{\bf c}$.

A regular subdivision of $A$ is called a {\em triangulation} if the
columns of each $A_{\sigma_i}$ are linearly independent for all
$i=1,\ldots,t$. Note that a regular subdivision $\Delta_{\bf c}$ is
a triangulation if and only if every vertex is contained in exactly
$d$ facets; that is, the polyhedron $Q_{\bf c}$ is {\em simple}, or,
{\em non-degenerate}. A triangulation $\Delta_{\bf c}$ is called
{\em unimodular} if $\textup{det}(A_{\sigma_i}) = \pm 1$ for each
maximal cell of $\Delta_{\bf c}$. A {\em refinement} of a
subdivision $\Delta_{\bf c}$ of $A$ is another subdivision
$\Delta_{{\bf c}^\prime}$ of $A$ so that each cell of $\Delta_{{\bf
c}^\prime}$ is contained in some cell of $\Delta_{\bf c}$. A set $B
\subset \Z^d$ is a {\em Hilbert basis} if $\N B = \textup{cone}(B)
\cap \Z^d$; we will say that a matrix $B$ is a Hilbert basis if its
columns form a Hilbert basis.  Note that if for some ${\bf c} \in
\R^n$ $\Delta_{\bf c}$ is a unimodular triangulation of $A$ then
Cramer's rule implies that $A$ itself is a Hilbert basis.

Let $ \textup{IP}_{A,{\bf c}} \, := \, \{ \textup{IP}_{A,{\bf
c}}({\bf b}) \, : \, {\bf b} \in \N A \} $ denote the family of
integer programs $\textup{IP}_{A,{\bf c}}({\bf b})$ having a
feasible solution. Informally, a {\em test-set} for the family of
integer programs $\textup{IP}_{A,{\bf c}}$ is a finite collection of
integer vectors, called {\em test vectors}, with the property that
any non-optimal feasible solution to any integer program in this
class can be improved (in objective value) by subtracting a test
vector from it. Test-sets for the family $\textup{IP}_{A,{\bf c}}$
were first introduced by Graver \cite{graver}. Graver's test set,
called a {\em Graver basis} in the literature, is typically not a
minimal test set.

A simple but helpful characterization of the TDI property in terms
of the Hilbert basis property of regular subdivisions has been
provided by Schrijver \cite{schrijver}. We prove another elementary
characterization in Section~\ref{sec:TDI} whose simplified
version is the following:

\medskip

\noindent {\bf Theorem~\ref{thm:equivalences}} {\em The system
${\bf y}A \leq {\bf c}$ is TDI if and only if
$A$ is a Hilbert basis, and there exists a test-set for
$\textup{IP}_{A,{\bf c}}$ with all test vectors having positive
entries equal to $1$ and this positive support indexes either
a linearly independent set or a minimally linear dependent set
of columns of $A$.}

\medskip

Establishing this theorem, some of its corollaries (see
Section~\ref{sec:TDI}) have led us to three relevant earlier results
that have been found in contexts different from ours, independently
of one another:

\begin{itemize}
\item[--]
Applegate, Cook \& McCormick's result \cite[Theorem 2]{ACM} which
states that if $A$ is a Hilbert basis then ${\bf y}A \leq {\bf c}$
is TDI if $\textup{LP}_{A,{\bf c}}({\bf b})$ has an integer optimum
for every constraint vector ${\bf b}$ that is a sum of linearly
independent columns of $A$. However, our test-sets and their use for
optimizing on IP are new.
\item[--]
Ho{\c s}ten and Sturmfels result \cite[Theorem 1.1]{gap} is stated
in a computer algebra context, but the special case when the gap is
$0$ can be seen to be equivalent to the above statement of
Theorem~\ref{thm:equivalences}. However, it is still easier to prove
each of these two independent results in their own terms, directly.
\item[--] Eisenbrand and Shmonin \cite{es} proved Ho{\c s}ten and
Sturmfels result in an elementary way, yet in the context of a
generalization of integer programming called ``parametric integer
programming''. However, this work also uses more involved machinery
than we do here.
\end{itemize}

In addition, none of these works characterizes the TDI property in terms
of test-sets or solving the underlying integer programs in fixed
dimension without using Lenstra's algorithm.
As previously noted in \cite{ACM}, this allows us to
deduce shortly Cook, Lov\'asz \& Schrijver's result \cite{cls} on
testing for the TDI property in fixed dimension, when given the
conditions of Theorem~\ref{thm:equivalences}.

Another virtue of Theorem~\ref{thm:equivalences} is a useful
reformulation to polynomial ideals: it generalizes a
well-known algebraic result proved by Sturmfels \cite[Corollary 8.9]{stu} 
relating toric initial ideals to unimodular 
triangulations. The basic connections between integer programming
and computational algebra, knowledge of which will not be assumed
here,  was initiated by Conti and Traverso
\cite{contitraverso} and further studied from various viewpoints in
\cite{swz}, \cite{stu}, \cite{recipes}, \cite{variations}. 

If $A$ is a matrix whose first $d \times (n-d)$ submatrix is a $0-1$
matrix and whose last $d\times d$ submatrix is $-I_d$, and ${\bf c}$
is all $1$ except for the last $d$ coordinates which are $0$, then
$\textup{DP}_{A,{\bf c}}({\bf b})$ is called a {\em set packing
problem}, and $Q_{\bf c}$ a {\em set packing polytope}.
In Section~\ref{sec:perfect} we show that the converse of the
following fact (explained at the end of Section~\ref{sec:TDI})
holds for normal fans of integral set packing polytopes:
{\em if
${\bf c}, {\bf c}^\prime \in \R^n$ are such that $\Delta_{{\bf
c}^\prime}$ is a refinement of $\Delta_{\bf c}$, where $\Delta_{{\bf
c}^\prime}$ is  a unimodular triangulation, then ${\bf y}A \leq {\bf
c}$ is TDI.} In general, the converse does not hold and the
most that is known in this direction is the existence
of just one full dimensional subset of the columns of $A$ which is
unimodular \cite{gs}. Not even a ``unimodular covering'' of a
Hilbert basis may be possible \cite{BRUNS}. However, the converse
does hold for normal fans of integral set packing polytopes. More
precisely, the main result of Section~\ref{sec:perfect} is the
following:

\medskip

\noindent {\bf Theorem~\ref{thm:perfect}} {\sl Given a set-packing
problem defined by $A$ and ${\bf c}$, $Q_{\bf c}$ has integer
vertices if and only if there exists ${\bf c}^\prime$ such that
$\Delta_{{\bf c}^\prime}$ is a refinement of the normal fan
$\Delta_{{\bf c}}$ of $Q_{\bf c}$, where $\Delta_{{\bf c}^\prime}$
is a unimodular triangulation.}

\medskip

The proof relies on the basic idea of Fulkerson's famous
``pluperfect graph theorem'' \cite{fulkerson} stating that the
integrality of such polyhedra implies their total dual integrality
in a very simple ``greedy'' way. Chandrasekaran and Tamir \cite{ct}
and Cook, Fonlupt and Schrijver \cite{cfs} exploited Fulkerson's
method by pointing out its lexicographic or advantageous
Carath{\' e}odory feature. In \cite[\S 4]{sebo} it is noticed with the
same method that the active rows of the dual of integral set packing
polyhedra (the cells of their normal fan) have a unimodular
subdivision, which can be rephrased as follows: {\em  the normal fan
of integral set packing polyhedra has a unimodular refinement.}
However, the proof of the regularity of such a refinement appears
for the first time in the present work.

These results offer four methods for recognizing TDI systems,
explained, illustrated and compared to previously known ones as
well, in Section~\ref{sec:computation}. Particular attention will be
given to the recognition of integral set packing polytopes in
practice, through Theorem~\ref{thm:perfect} above, and by using
computational algebra packages like {\tt Macaulay 2}.

\section{TDI Systems} \label{sec:TDI}

In this section we provide some new characterizations of TDI
systems. We show the equivalence of five properties, three
polyhedral (one of them is the TDI property) and two concern
polynomial ideals. A third property is also equivalent to these in
the generic case.

While the proofs of the equivalences of the three polyhedral
properties use merely polyhedral arguments, the last among them
-- (iii) -- has an appealing reformulation into the
language of polynomial ideals.
Therefore, we start this section by introducing the necessary
background on polynomial ideals; namely, toric ideals, their
initial ideals and Gr{\" o}bner bases. The characterizations
of TDI systems involving polynomial ideals are useful generalizations
of known results in computational algebra.
See \cite{clos1} and \cite{stu} for further background.

An {\em ideal} $I$ in a polynomial ring $R := {\bf k}[x_1, \ldots,
x_n]$ is an $R$-vector subspace with the property that $I \cdot R =
I$. It was proven by Hilbert that every ideal is finitely generated.
That is, given an ideal $I$ there exists a finite set of polynomials
$f_1, \ldots, f_t \in I$ such that for every $f \in I$ there exists
$h_1, \ldots, h_t \in R$ with $f = h_1 f_1 + \cdots + h_t f_t$. We
call such a collection $f_1, \ldots, f_t \in I$ a {\em generating
set} for the ideal $I$ and denote this by $ I = \langle f_1, \ldots,
f_t \rangle. $ For the monomials in $R$ we write ${\bf x}^{\bf u} =
x_1^{u_1} \cdots x_n^{u_n}$ for the sake of brevity. We call ${\bf
u}$ the {\em exponent vector} of ${\bf x}^{\bf u}$. A monomial ${\bf
x}^{\bf u}$ is said to be {\em square-free} if ${\bf u} \in
\{0,1\}^n$. An ideal is called a {\em monomial ideal} if it has a
generating set consisting only of monomials. For any ideal $J$ of
$R$, $\textup{mono}(J)$ denotes the monomial ideal in $R$ generated
by the set of monomials in $J$. An algorithm  for computing the
generators of the monomial ideal $\textup{mono}(J)$ can be found
in \cite[Algorithm 4.4.2]{sst}.

Every weight vector ${\bf c} \in \R^n$ induces a
partial order $\succeq$ on the monomials in $R$ via ${\bf x}^{\bf u}
\succeq {\bf x}^{\bf v}$ if ${\bf c} \cdot {\bf u} \geq {\bf c}
\cdot {\bf v}$. If ${\bf c} \in \R^n$ where $1$ is the monomial of
minimum ${\bf c}$-cost (that is, ${\bf c}\cdot {\bf u} \ge 0$ for
every monomial ${\bf x}^{\bf u})$, then we can define initial
terms and initial ideals. Given a
polynomial
$f = \sum_{{\bf u} \in \N^n} r_{\bf u}{\bf x}^{\bf u} \in I$
the {\em initial term} of $f$ with respect to ${\bf c}$, is
denoted by $ \textup{in}_{\bf c}(f) $, and equals the sum of all
$r_{\bf u}{\bf x}^{\bf u}$ of $f$, where ${\bf c} \cdot {\bf u}$
is maximum.
Note that we can always write a polynomial $f$ as
$f = \textup{in}_{\bf c}(f) + \textup{trail}_{\bf c}(f)$
in the obvious way.
The {\em initial ideal} of $I$
with respect to ${\bf c}$ is defined as the ideal in $R$ generated
by the initial terms of the polynomials in $I$: $ \textup{in}_{\bf
c}(I) := \langle \, \textup{in}_{\bf c}(f) \, : \, f \in I \,
\rangle$. A {\em Gr{\" o}bner basis} of an ideal $I$ {\em with respect to
${\bf c}$}, is a finite collection of elements $g_1, \ldots, g_s$ in
$I$ such that
$ \textup{in}_{\bf c}(I) = \langle \,
\textup{in}_{\bf c}(g_1), \, \textup{in}_{\bf c}(g_2), \, \ldots,
\textup{in}_{\bf c}(g_s) \, \rangle$. Every Gr{\" o}bner basis
is a generating set for the ideal $I$.

If $\textup{in}_{\bf c}(I)$ is a monomial ideal then a
Gr{\" o}bner basis is {\em reduced} if for every $i \neq j$,
no term of $g_i$ is divisible by $\textup{in}_{\bf c}(g_j)$.
The reduced Gr{\" o}bner basis is unique. In this case,
the set of monomials in $\textup{in}_{\bf c}(I)$ equal
$\{{\bf x}^{\bf u}: {\bf u} \in U \}$ with $U:= D + \N^n$ where
$D$ is the set of exponent vectors of the monomials
$\textup{in}_{\bf c}(g_1), \, \textup{in}_{\bf c}(g_2), \, \ldots,
\textup{in}_{\bf c}(g_s)$. Dickson's lemma states that sets of the
form $D + \N^n$, where $D$ is
arbitrary have only a finite number of minimal elements (with
respect to coordinate wise inequalities). This is an alternative proof
to Hilbert's result that every polynomial ideal is finitely generated.
In this case, the Gr{\" o}bner basis also provides a generalization
of the Euclidean algorithm for polynomial rings with two or more variables
called Buchberger's algorithm (see \cite[Chapter 2, \S 7]{clos1}).
This algorithm solves
the {\em ideal membership problem}: decide if a given polynomial is
in an ideal or not.
However, a Gr{\" o}bner basis for an ideal can have many
elements (compared to a minimal generating set for the ideal), and
none of the related computations can be achieved in polynomial time.

If $\textup{in}_{\bf c}(I)$ is not a monomial ideal then we can form
a monomial initial ideal, $\textup{in}_{\succ_{\bf c}}(I)$, as
follows: fix an arbitrary {\em term order} independent of ${\bf c}$,
(that is, a total ordering $\succ$  of the vectors in $\N^n$
satisfying ${\bf u} \succ {\bf 0}$ for every ${\bf u} \in \N^n$, and
if ${\bf u} \succ {\bf v}$ then ${\bf u} + \gamma \succ {\bf v} +
\gamma $ for all $\gamma \in \N^n$). We use this term order to break
ties: ${\bf x}^{\bf u} \succ_{\bf c} {\bf x}^{\bf v}$ if and only if
${\bf c} \cdot {\bf u} > {\bf c} \cdot {\bf v}$, or ${\bf c} \cdot
{\bf u} = {\bf c} \cdot {\bf v}$ and ${\bf u} \succ {\bf v}$.
Clearly, ``$\succ_{\bf c}$'' is also a term order and so
$\textup{in}_{\succ_{\bf c}}(I)$ is a monomial ideal. Such a
tie-breaking will be needed in the proof of Lemma~\ref{lem:mono}.

The {\em toric ideal} of $A$ is the ideal $ I_A = \langle {\bf
x}^{\bf u} - {\bf x}^{\bf v} \, : \, A{\bf u} = A{\bf v}, \, \, {\bf
u}, {\bf v} \in \N^n \, \rangle $ and is called a {\em binomial
ideal} since it is generated by polynomials having at most terms.
Every reduced Gr{\" o}bner basis of a toric ideal consists of
binomials. A {\em toric initial ideal} is any initial ideal of a
toric ideal.

\begin{remark} \label{rem:gb}
It follows for (a not necessarily generic) ${\bf c}$, that the
reduced Gr\"obner basis of a toric ideal $I_A$ (with respect to
$\succ_{\bf c}$) is of the form $\G_{\succ_{\bf c}} = \{ {\bf
x}^{{\bf u}_i^+} - {\bf x}^{{\bf u}_i^-} \, : \, i = 1, \ldots, t
\}$ and we can suppose that $ \textup{in}_{\bf c}({\bf x}^{{\bf
u}_i^+} - {\bf x}^{{\bf u}_i^-}) = {\bf x}^{{\bf u}_i^+} $ for $i
=1, \dots, s$ and $ \textup{in}_{\bf c}({\bf x}^{{\bf u}_i^+} - {\bf
x}^{{\bf u}_i^-}) = {\bf x}^{{\bf u}_i^+} - {\bf x}^{{\bf u}_i^-} $
for $i =s+1, \dots, t$. Furthermore,  the set of polynomials $
\Set_{\succ} :=\{ {\bf x}^{{\bf u}_1^+}, \ldots, {\bf x}^{{\bf
u}_s^+}, {\bf x}^{{\bf u}_{s+1}^+} - {\bf x}^{{\bf u}_{s+1}^-},
\ldots, {\bf x}^{{\bf u}_t^+} - {\bf x}^{{\bf u}_t^-} \} $ is a
reduced Gr{\" o}bner basis for $\textup{in}_{\bf c}(I_A)$ with
respect to the term order $\succ$  cf. \cite[Corollary 1.9]{stu}.
Note that if ${\bf c}$ were generic then $s=t$ and 
$\Set_\succ$ is simply the minimal generating set for
$\textup{in}_{\bf c}(I_A)$.
\end{remark}

The following lemma is a natural connection between integer programming
and toric initial ideals. It originally appeared in \cite[Lemma 4.4.7]{sst}
but we prove it here in order for this article to be self-contained.

\begin{lemma} \label{lem:mono}
For $A \in \Z^{d \times n}$ and ${\bf c} \in \R^n$ the monomial
ideal $\textup{mono}(\textup{in}_{\bf c}(I_A))$ is equal to
$$\langle \, {\bf x}^\omega \, : \, \omega\in\N^n \, \textup{is a non-optimal
solution for} \, \textup{IP}_{A,{\bf c}}(A \omega) \, \rangle.
$$
\end{lemma}

\begin{proof}
It is straightforward to show that $\textup{mono}(\textup{in}_{\bf c}(I_A))$
contains the defined set: let $\omega$  be a non-optimal solution,
and $\omega'$ an optimal solution to $\textup{IP}_{A,{\bf c}}(A
\omega)$. Then ${\bf x}^\omega - {\bf x}^{\omega'} \in  I_A$ is a
binomial having the monomial ${\bf x}^\omega\in
\textup{mono}(\textup{in}_{\bf c}(I_A))$ as its initial term with
respect to ${\bf c}$.

Suppose now ${\bf x}^\omega \in \textup{in}_{\bf c}(I_A)$.
That is, there exists a polynomial $f \in I_A$ with
$\textup{in}_{\bf c}(f) = {\bf x}^\omega$. We will prove that
$\omega$ is a non-optimal solution for
$\textup{IP}_{A,{\bf c}}(A \omega)$.

Let $\succ$ be an arbitrary term order and let $\Set_{\succ}$
be the reduced Gr{\" o}bner basis of $\textup{in}_{\bf c}(I_A)$
with respect to $\succ$ as in Remark~\ref{rem:gb}. We proceed by induction with
respect to the minimum number of successive polynomial divisions
of $f$ by the elements of $\Set_{\succ}$, and replacing
$f$ by the remainder after each division, until arriving at a $0$ remainder.
As noted in Remark~\ref{rem:gb}
$\Set_{\succ} =\{
{\bf x}^{{\bf u}_1^+}, \ldots, {\bf x}^{{\bf u}_s^+}, {\bf x}^{{\bf
u}_{s+1}^+} - {\bf x}^{{\bf u}_{s+1}^-}, \ldots, {\bf x}^{{\bf
u}_t^+} - {\bf x}^{{\bf u}_t^-} \} $.


If first we divide $f$ with a monomial in $\Set_\succ$, say ${\bf x}^{{\bf
u}_1^+}$, then ${\bf c} \cdot {\bf u}_1^+>{\bf c} \cdot {\bf
u}_1^-$, and $\textup{in}_{\bf c}(f) = {\bf x}^\omega$ is divisible
by ${\bf x}^{{\bf u}_1^+}$, so $\omega - {\bf u}_1^+ \geq {\bf 0}$,
$ {\bf c} \cdot (\omega - ({\bf u}_1^+ - {\bf u}_1^-)) < {\bf c}
\cdot \omega $, where $A(\omega - ({\bf u}_1^+ - {\bf
u}_1^-))=A\omega$. Threrefore $\omega$ is also a non-optimal
solution to $ \textup{IP}_{A, {\bf c}}(A\omega) $.

Otherwise, the first polynomial division is by a binomial in $\Set_\succ$,
say, ${\bf x}^{{\bf u}_t^+} - {\bf x}^{{\bf u}_t^-}$.
Now recall that $f \in I_A$ is a polynomial with
$\textup{in}_{\bf c}(f) = \textup{in}_{\succ_{\bf c}}(f)$ and so the resulting
polynomial after the division can be written as
$ f^\prime :=
\frac{{\bf x}^\omega}{{\bf x}^{{\bf u}_t^+}} {\bf x}^{{\bf u}_t^-}
\, - \, \textup{trail}_{{\bf c}}(f) $.
The initial term $\textup{in}_{\bf c}(f^\prime)$ is a monomial and equals
$
\frac{{\bf x}^\omega}{{\bf x}^{{\bf u}_t^+}} {\bf x}^{{\bf u}_t^-}
$
since ${\bf c}\cdot(\omega - ({\bf u}_t^+ - {\bf u}_t^-))$ equals
${\bf c} \cdot \omega$, due to ${\bf c} \cdot {\bf u}_t^+ = {\bf c} \cdot {\bf u}_t^-$.

This implies that
$\textup{in}_{\bf c} (\textup{trail}_{{\bf c}}(f))$
is strictly cheaper (with respect to ${\bf c}$) than
$
\frac{{\bf x}^\omega}{{\bf x}^{{\bf u}_t^+}} {\bf x}^{{\bf u}_t^-}
$
and so
$
\frac{{\bf x}^\omega}{{\bf x}^{{\bf u}_t^+}}
{\bf x}^{{\bf u}_t^-} \in \textup{in}_{\bf c}(I_A)
$.
Since $f^\prime$ requires one less division than $f$ by
$\Set_\succ$ to arrive at $0$ remainder then, by induction,
we have that $\omega - ({\bf u}_t^+ - {\bf u}_t^-)$ is not an optimal
solution to
$\textup{IP}_{A, {\bf c}}(A(\omega - ({\bf u}_t^+ - {\bf u}_t^-)))$
and so $\omega$ is not an optimal solution to
$\textup{IP}_{A, {\bf c}}(A\omega)$ either.
\end{proof}

\begin{remark} \label{rem:supp}
$\sigma \subseteq [n]$ is contained in a cell of $\Delta_{\bf c}$
if and only if $\sum_{i \in \sigma} {\bf e}_i \in \R^n$ is an
optimal solution to $ \textup{LP}_{A, {\bf c}}({\bf b}_{\sigma}) $
where ${\bf b}_{\sigma} := A(\sum_{i \in \sigma} {\bf e}_i)$.
This happens in particular if $\sigma=\emptyset$. Indeed, then
$\sum_{i \in \sigma} {\bf e}_i={\bf 0}$, and
${\bf b}_{\sigma} := A(\sum_{i \in \sigma} {\bf e}_i)={\bf 0}$.
For this ${\bf b}={\bf 0}\in\R^d$, ${\bf 0}\in \R^n$ is an optimal
solution, otherwise $\textup{LP}_{A,{\bf c}}({\bf b})$ is not bounded
for any ${\bf b}$, and $Q_c=\emptyset$. We are not interested in such
polyhedra, and we will avoid this situation by supposing
$Q_{\bf c} \ne \emptyset$.
\end{remark}

A {\em test-set} \cite{graver} for the family of integer programs
$\textup{IP}_{A,{\bf c}} \, := \,
\{ \textup{IP}_{A,{\bf c}}({\bf b}) \, : \, {\bf b} \in \N A \}
$
is a collection of integer vectors
$
\{ {\bf v}_i^+ - {\bf v}_i^- \, : \,
A{\bf v}_i^+ = A{\bf v}_i^-, \,
{\bf c\cdot v_i}^+ > {\bf c\cdot v_i}^-,\,\,
{\bf v}_i^+, {\bf v}_i^- \in \N^n,\, i = 1, \ldots, s \}
$
with the property that for all feasible, non-optimal solution
${\bf u}$ to $\textup{IP}_{A,{\bf c}}({\bf b})$  there exists an
$i$, $1 \leq i \leq s$, such that ${\bf u} - ({\bf v}_i^+ - {\bf
v}_i^-) \geq {\bf 0}$. We can now state and prove our
characterizations of TDI.

\begin{example} \label{ex:test}
Let us suppose $Q_{\bf c} \ne \emptyset$, that is, ${\bf 0}\in \R^n$
is an optimal solution for $\textup{LP}_{A,{\bf c}}({\bf 0})$, and
show a particular test-set related to Remark~\ref{rem:supp}. We will
say that $\kappa \subseteq [n]$ is a {\em wheel} (with respect to
$A$ and ${\bf c}$), if $A(\sum_{i \in \kappa} {\bf e}_i)={\bf 0}$,
$\sum_{i \in \kappa} {\bf c}_i > 0$ and for all $i \in \kappa$,
$\kappa \setminus\{i\}$ is a subset of a cell (of
$\Delta_{\bf c}(A)$).

For instance if  $A=({\bf a}_1,\ldots, {\bf a}_d, -({\bf
a}_1+\ldots+{\bf a}_d))$, where ${\bf a}_1,\ldots, {\bf a}_d$ are
linearly independent integer vectors, and ${\bf c} \ge 0$, then
$\Delta_{\bf c}(A)$ is a triangulation whose maximal cells are
precisely $\{ [d+1] \backslash \{ i \} \, : \, i = 1,2 \ldots d+1
\}$. In this case the one-element set $\{(1,\ldots,1)\in\R^{d+1}\}$
is a test-set for $\textup{IP}_{A, {\bf c}}$. Note that if
$\textup{cone}(A)$ is pointed then $A$ has no wheel, regardless of
the ${\bf c}$ in question.
\end{example}

\begin{theorem} \label{thm:equivalences}
Fix $A \in \Z^{d \times n}$, where $A$ is a Hilbert basis,
${\bf c} \in \R^n$, and $Q_{\bf c}\ne\emptyset$.
The following statements are equivalent:
\begin{enumerate}
\item[(i)] The system ${\bf y}A \leq {\bf c}$ is TDI.
\item[(ii)] The subconfiguration $A_\sigma$ of $A$ is a Hilbert basis
for every cell $\sigma$ in $\Delta_{\bf c}$.
\item[(iii)] There exists a test-set for $\textup{IP}_{A,{\bf c}}$
where all the positive coordinates are equal to $1$, and each
positive support is either the incidence vector of a linearly
independent set of columns, or  of a wheel; in the former case the
negative support is a subset of a cell, in the latter case it is
$\emptyset$.
\item[(iv)] The monomial ideal
$
\langle \, {\bf x}^\omega : \omega \in \N^n \, \textup{is not an optimal
solution for} \, \textup{IP}_{A,{\bf c}}(A \omega) \, \rangle
$
has a square-free generating set.
\item[(v)] The monomial ideal generated by the set of monomials
in  $\textup{in}_{\bf c}(I_A)$,
that is,  $\textup{mono}(\textup{in}_{\bf c}(I_A))$ has a
square-free generating set.
\end{enumerate}
\end{theorem}

The main content of the equivalence of (ii) and (v) is that the
``Hilbert basis property'' is equivalent to the existence of a
``square-free generating set'', extending a well-known result
\cite[Corollary 8.9]{stu} (see Corollary~\ref{thm:sturm} below) to
the non-generic case. We have recently found a similar extension in
Ho{\c s}ten and Sturmfels \cite{gap}, equivalent to ours. 
(The latter paper provides 
 an algorithm to compute the { integer programming gap} $ \textup{gap}_{A,{\bf c}}
:= \textup{max} \{ \textup{OPTIP}_{A,{\bf c}}({\bf b}) -
\textup{OPTLP}_{A,{\bf c}}({\bf b}) \, : \, {\bf b} \in \N A \} $,
where OPTIP and OPTLP mean the optimal value of the corresponding programs.
Note that the system ${\bf y}A \leq {\bf c}$ being TDI is equivalent
to $A$ being a Hilbert basis and $\textup{gap}_{A,{\bf c}} = 0$.
Ho{\c s}ten and Sturmfels compute $\textup{gap}_{A,{\bf c}}$ by
studying the ``primary decomposition'' of the monomial ideal
$\textup{mono}(\textup{in}_{\bf c}(I_A))$. We present below our
direct elementary argument.)

Note also that (i) implies that $A$ is a Hilbert
basis, since this latter is equivalent to the following property: if
$\textup{LP}_{A,{\bf c}}({\bf b})$ has a feasible solution, then it also
has an integer feasible solution; (ii) also implies it, since a cone
that has a subdivision to Hilbert cones is itself a Hilbert basis.
So this condition has to be added only to (iii), (iv) and (v); they
may all be satisfied without $A$ being a Hilbert basis, certifying
then that ${\bf y}A \leq {\bf c}$ is not TDI, unless $Q_{\bf c}=\emptyset$.
The condition $Q_{\bf c}\ne\emptyset$ makes sure that
$\textup{LP}_{A,{\bf c}}({\bf b})$ is bounded, otherwise already the
statement of (iii), (iv), (v) is not well defined.

\begin{proof}
{\bf (i) $\Rightarrow$ (ii) }: This is well-known from Schrijver's
work, (see for instance \cite{schrijver}), but we provide the (very simple)
proof here for the sake of completeness: Suppose the system
${\bf y}A \leq {\bf c}$ is TDI, and let $\sigma\in\Delta_c$.
We show that $A_\sigma$ is a Hilbert basis.
Let ${\bf b} \in \textup{cone}(A_{\sigma})$. Since the optimal solutions
for $\textup{LP}_{A, {\bf c}}({\bf b})$ are exactly the non-negative
combinations  of the columns of $\A_{\sigma}$ with result ${\bf b}$,
the TDI property means exactly that ${\bf b}$ can also be written as a
non-negative integer combination of columns in $A_{\sigma}$, as claimed.

\vspace{.2cm}

\noindent {\bf (ii) $\Rightarrow$ (iii) }: Suppose (ii) holds true
for $\Delta_{\bf c}$ of $A$. For every $\tau \subseteq [n]$ with
$\tau$ not contained in any cell of $\Delta_{\bf c}$, let ${\bf
b}_\tau :=\sum_{i \in \tau}{\bf a}_i = A(\sum_{i \in \tau} {\bf
e}_i)$. Since $\tau$ is not contained in any cell of $\Delta_{\bf
c}$, there exists an optimal solution $\beta_\tau$ to
$\textup{LP}_{A, {\bf c}}({\bf b}_{\tau}) $ with ${\bf c} \cdot
\beta_\tau < {\bf c} \cdot \sum_{i \in \tau} {\bf e}_i$. By the
optimality of $\beta_\tau$ we must have $\textup{supp}(\beta_\tau)
\subseteq \sigma$ for some cell $\sigma$ of $\Delta_{\bf c}(A)$.

If ${\bf b}_\tau \ne {\bf 0}$, then $\sigma\ne\emptyset$, and by (ii)
$A_\sigma$ is a Hilbert basis. Therefore $\beta_\tau$ can be chosen
to be an integral vector. On the other hand, if ${\bf b}_\tau = {\bf 0}$,
then by the condition $Q_{\bf c}\ne\emptyset$, ${\bf 0}$ is an optimal
solution for $\textup{LP}_{A, {\bf c}}({\bf b}_{\tau}) $
(see Remark~\ref{rem:supp}). Let
\begin{align*}
{\mathcal T}_{A,{\bf c}} \, := & \, \{ \,
\sum_{i \in \tau} {\bf e}_i - \beta_\tau \, : \,
A_\tau \, \textup{is linearly independent and}
 \, \tau \, \textup{is not contained in any cell of} \, \Delta_{\bf c} \, \} \\
& \, \bigcup \, \{\sum_{i \in \kappa} {\bf e}_i : \, \kappa \,
\textup{ is a wheel (with respect to $A$ and ${\bf c}$)} \}.
\end{align*}

We claim that ${\mathcal T}_{A,{\bf c}}$ is a test-set
for $\textup{IP}_{A, {\bf c}}$. By construction every
${\bf t} \in {\mathcal T}_{A,{\bf c}}$ satisfies $A{\bf t} = {\bf 0}$
and ${\bf c} \cdot {\bf t} >0$, so we have to prove only the following:

\vspace{.1cm}

\noindent{\bf Claim}: {\em For every ${\bf b} \in \Z^d$ and feasible
but not optimal solution $\omega$ of $\textup{IP}_{A, {\bf c}}({\bf
b})$, there exists ${\bf t} \in {\mathcal T}_{A,{\bf c}}$ such that
$\omega - {\bf t} \ge {\bf 0}$.}

We have $A\omega = {\bf b}$, where $\textup{supp}(\omega)$ is not
contained in any cell $\sigma$ of $\Delta_{\bf c}$. 

\vspace{.1cm}

\noindent{\bf Case 1.} ${\bf b}:=A\omega\ne {\bf 0}$:

By basic linear programming (``Caratheodory's theorem'') there
exists $\tau \subseteq \textup{supp}(\omega)$, so that $A_\tau$ is
linearly independent, and ${\bf b}\in\textup{cone}(A_\tau)$. If
every  $\tau$ satisfying these three conditions is a subset of a
cell, then the unique solution of the equation $A_\tau {\bf x}_\tau
= {\bf b}$ is feasible for $\textup{LP}_{A, {\bf c}}({\bf b})$ and
optimal, and again, by basic linear programming this would
contradict that $\omega$ is not optimal. So $\tau$ can be chosen so
that it satisfies these three properties, but {\em it is not
contained in any cell of} $\Delta_{\bf c}$. Clearly, $ \omega -
(\sum_{i \in \tau} {\bf e}_i - \beta_\tau) \geq {\bf 0}$, and ${\bf
t} := \sum_{i \in \tau} {\bf e}_i - \beta_\tau\in {\mathcal
T}_{A,{\bf c}}$, finishing the proof of the claim.

\noindent{\bf Case 2.} ${\bf b}:=A\omega= {\bf 0}$:

\vskip 0.1cm Since then $\bf 0$ is optimal, and $\omega$ is not, we
have $${\bf c} \cdot \omega >0. \eqno{(1)}$$

\noindent (a) {\em every proper subset of $\textup{supp}(\omega)$ is
contained in some cell, and $\omega$ is a $0-1$ vector.}

\vspace{.1cm}

First, if there exists $\omega' \in \N^n$, $\omega' \le \omega$,
with $\omega'$ not optimal for ${\bf 0} \neq {\bf b}':= A \omega'$,
then we can apply Case 1 to $\omega'$, and get $\bf t\in {\mathcal
T}_{A,{\bf c}}$, $\omega-t\ge \omega'-t\ge 0$, and we are done.
Second, if $\omega$ has a coordinate that is bigger than $1$, say
$\omega_1
>1$, then Case~1 can be applied to $\omega':=\omega - e_1$.

We can now assume that (a) holds. Letting
$\kappa:=\textup{supp}(\omega)$, we show:

\vspace{.1cm}

\noindent (b) {\em every proper subset of $A_\kappa$ is linearly
independent.}

\vspace{.1cm}

Suppose (b) is not true. Then for some $J\subsetneq\kappa$ there 
exists a coefficient vector 
$\lambda = (\lambda_1,\ldots,\lambda_n)\in\mathbb{R}^n$, 
$\lambda_j \ne 0$ if $j\in J$, and $\lambda_j= 0$ if $j\notin J$.  
Moreover (by dividing with the highest coefficient, say that of the first 
coordinate), $\lambda_j \le 1$ $(j\in [n])$, and $\lambda_1=1$, such 
that 
$$
A{\bf \lambda }={\bf 0}.
$$ 
Furthermore, we can partition $J = J^+ \cup J^-$ where 
$J^+ := \{ j : \lambda_j > 0  \}$ and 
$J^- := \{ j : \lambda_j < 0  \}$. 
Using this partition, we can rewrite 
$\lambda$ in the obvious way as $\lambda = \lambda^+ - \lambda^-$ 
where $\lambda^+, \lambda^- \geq {\bf 0}$ with the supports being 
$J^+$ and $J^-$ respectively. Consequently, we have 
$$
A \lambda^+ = A \lambda^-
$$ 
and according to (a) there exist cells of $\Delta_{\bf c}(A)$ 
containing $J^+$ and $J^-$ so both $\lambda^+$ and $\lambda^-$ 
are optimal solutions to the same linear program implying that:
$$
{\bf c} \cdot \lambda^+ ={\bf c} \cdot \lambda^-. 
\eqno{(2)}
$$

By our choice of the size of the components in 
$\lambda = \lambda^+ - \lambda^-$, we have $\lambda_1 = 1$. 
In addition, $\lambda \leq \omega$ and $\omega_1 = 1$ and, 
since 
$
\textup{supp}(\omega - (\lambda^+ - \lambda^-)) 
\subseteq
\kappa \backslash \{ 1 \}
\subsetneq 
\kappa
$ 
then, by (a) and the fact that $A\lambda^+ = A\lambda^-$, 
$\omega - (\lambda^+ - \lambda^-)$ is an optimal 
solution to $\textup{LP}_{A,{\bf c}}({\bf b})$. 

Recall that ${\bf 0}$ is also an optimal solution to 
$\textup{LP}_{A,{\bf c}}({\bf b})$. 
Combining (1) and (2) above, we get that 
${\bf c} \cdot (\omega - (\lambda^+ - \lambda^-)) > 0 = {\bf c} \cdot {\bf 0}$ 
which contradicts the optimality of $\omega - (\lambda^+ - \lambda^-)$ 
and so no such linear dependent $J \subsetneq \kappa$ can exist. This proves (b).

Now by (a) and (b):  $\kappa$ is a wheel for $A$ and ${\bf c}$, and as
such $\sum_{i \in \kappa} {\bf e}_i\in {\mathcal T}_{A,{\bf c}}$,
finishing the proof of the claim and of (ii) $\Rightarrow$ (iii).

\vspace{.2cm}

\noindent {\bf (iii) $\Rightarrow$ (i)} :  Suppose (iii) holds,
and  ${\bf b} \in \textup{cone}(A)$. Then $\textup{LP}_{A,{\bf c}}({\bf b})$
is feasible, and since $Q_{\bf c}\ne\emptyset$ the optimum is bounded.
So by the condition that $A$ is a Hilbert basis,
$\textup{IP}_{A,{\bf c}}({\bf b})$ is also feasible and it also has a bounded minimum
attained by  $\omega \in \N^n$.
 Suppose for a contradiction that the optimal solution $\alpha/D$ to
 $\textup{LP}_{A,{\bf c}}({\bf b})$,
 $( \alpha \in \N^n$, $D$ is a positive integer)
 satisfies ${\bf c} \cdot \alpha /D < {\bf c} \cdot \omega$.
This also implies that $D \omega$ is not an optimal solution to
$\textup{IP}_{A,{\bf c}}(D{\bf b})$.

By (iii) there exists a $\gamma^+ - \gamma^-\in {\mathcal T}_{A,{\bf
c}}$ with $\gamma^+ \in \{0,1\}^n$ and $\gamma^- \in \N^n$ such that
${\bf c} \cdot (\gamma^+ - \gamma^-) > 0$ and  $D \omega - (\gamma^+
- \gamma^-) \in \N^n$. Hence, $\textup{supp}(\gamma^+) \subseteq
\textup{supp}(D \omega) = \textup{supp}(\omega)$. Since the value of
all elements in $\gamma^+$ is $0$ or $1$ then we also have $\omega
\geq \gamma^+$, so $\omega - (\gamma^+ - \gamma^-) \in \N^n$  is
also a feasible solution to $\textup{IP}_{A,{\bf c}}({\bf b})$ with
${\bf c} \cdot (\omega - (\gamma^+ - \gamma^-) ) < {\bf c} \cdot
\omega$, in contradiction to the optimality of $\omega$.

\vspace{.2cm}

\noindent {\bf (iii) $\Leftrightarrow$ (iv) $\Leftrightarrow$ (v):}
Both (iii) and (iv) can
be reformulated as follows: If $\omega \in \N^n$ is not an optimal
solution to $\textup{IP}_{A,{\bf c}}(A \omega)$ then the vector
$\omega^\prime := \sum_{i \in \textup{supp}(\omega)} {\bf e}_i$
is also a non-optimal solution to
$\textup{IP}_{A,{\bf c}}(A \omega^\prime)$. The equivalence of
(iv) and (v) is a special case of Lemma~\ref{lem:mono}.
\end{proof}

The simple equivalence of (i) with (iii) is a possible alternative
to Cook, Lov\'asz \& Schrijver's  result \cite{cls}, and (iii) could
also be replaced by a characterization of Applegate, Cook \&
McCormick \cite{ACM} stating that it is sufficient to check the
property for functions ${\bf b}$ that are small subsums of rows of
$A$. For the sake of completeness and for later reference we present
the corresponding computational corollary in our context:

\begin{corollary} \cite{ACM}, \cite{cls} \label{thm:cls}
Let the dimension $d$ be fixed but the system
${\bf y}A \leq {\bf c}$, given as input, have an arbitrary number of inequalities, where
$A$ is a Hilbert basis. Then this system can be tested for the TDI
property in polynomial time.
\end{corollary}

\begin{proof}  We use the equivalence of (i) and (iii)
and the construction of ${\mathcal T}_{A,{\bf c}}$ from above.
Note that the wheels are straightforward to identify but they
do not have to be identified for testing TDI-ness.

We can construct the non-wheels of ${\mathcal T}_{A,{\bf c}}$
or conclude that the system is not TDI in $O(n^d)$ time. Listing all the linearly
independent (or generously all the $d$ element) subsets $\tau$ of
$\{1,\ldots,n\}$ that are not subsets of cells in the fan.
For each such subset we can form the vector
${\bf b}_\tau:= A(\sum_{i \in \tau} {\bf e}_i)$ and identify its cell
$\sigma$ in the fan by solving $\textup{LP}_{A,{\bf c}}({\bf b}_\tau)$.

Next, for each such $\tau$, with its corresponding optimal cell
$\sigma$, do the following: we can either find an integer optimal
solution $\beta_\tau$ to $\textup{LP}_{A,{\bf c}}({\bf b}_\tau)$
(which will have support in $\sigma$), or else no such integer
solution exists in which case we conclude that ${\bf y}A \leq {\bf
c}$ is not TDI, which can be done in polynomial time -- see \cite[\S
6.7]{gls} for example.

If $\textup{LP}_{A,{\bf c}}({\bf b}_\tau)$ has an integer solution
then repeat the same test, as for ${\bf b}_\tau$ above, for ${\bf
b}_\tau - r_i{\bf a}_i$  where $i \in \sigma$ and $r_i$ is the
largest positive integer for which the optimal solution to $
\textup{LP}_{A,{\bf c}}({\bf b}_\tau - r_i{\bf a}_i) $ is in the
same cell as $\sigma$. If $ \textup{LP}_{A,{\bf c}}({\bf b}_\tau -
r_i{\bf a}_i) $ has an integer optimum solution then repeat as above
by replacing ${\bf b}_\tau$ with ${\bf b}_\tau - r_i{\bf a}_i$. If
at any point during this procedure the updated ${\bf b}_\tau$ yields
no integer optimal solution to $\textup{LP}_{A,{\bf c}}({\bf
b}_\tau)$ then we immediately stop and conclude that ${\bf y} A \leq
{\bf c}$ is not TDI. Otherwise, we add $\sum_{i \in \tau} {\bf e}_i
- \beta_\tau$ to the test set.

We repeat this procedure for every such $\tau$ as above and
if for none of the linearly independent sets of columns $\tau$
(different from the cells) we arrive at the conclusion that the
system is not TDI, then we have in ${\mathcal T}_{A,{\bf c}}$ an
element of the form $\sum_{i \in \tau} {\bf e}_i - \beta_\tau$ for
every $\tau$ that is not subset of a cell, and we can conclude like
in the proof of the theorem, that ${\mathcal T}_{A,{\bf c}}$ is a
test-set. Therefore (iii) is satisfied, so according to
Theorem~\ref{thm:equivalences}, (i) is also satisfied, that is,
${\bf y}A \leq {\bf c}$ is TDI. If $d$ is fixed, $A$ has a
polynomial number of linearly independent subsets of columns.
\end{proof}

One consequence of this corollary is that the constructed test-set
allows us to solve $\textup{IP}_{A,{\bf c}}({\bf b}) $ in polynomial
time without the use of Lenstra's algorithm.

Recall that we defined  ${\bf c} \in \R^n$ to be {\em generic}
with the first of the following conditions; the others are
equivalent properties for toric ideals \cite{variations}:
\begin{enumerate}
\item[--] The integer program
$\textup{IP}_{A, {\bf c}}({\bf b})$ has a unique optimal solution
for all ${\bf b} \in \N A$.
\item[--] The toric initial ideal $\textup{in}_{\bf c}(I_A)$ is a monomial ideal.
\item[--] There exists a reduced Gr{\" o}bner basis
$
\{
{\bf x}^{{\bf u}_1^+} - {\bf x}^{{\bf u}_1^-},
\ldots,
{\bf x}^{{\bf u}_s^+} - {\bf x}^{{\bf u}_s^-}
\}
$
of $I_A$ with
$
{\bf c} \cdot {\bf u}_i^+ > {\bf c} \cdot {\bf u}_i^-
$
for each $i = 1,\ldots, s$.
\end{enumerate}

Recall that if each $\textup{IP}_{A, {\bf c}}({\bf b})$ has a
unique solution then so does $\textup{LP}_{A, {\bf c}}({\bf b})$.
Hence in the generic case, assuming $Q_{\bf c} \ne \emptyset$,
Cramer's rule tells us that $\Delta_{\bf c}$ being a unimodular
triangulation of $A$ implies that $A$ is a Hilbert basis
{\em without supposing it in advance}; and so (ii) implies (v)
without any condition; to prove the converse
statement in the generic case, (v) along with the assumption that
$\Z A = \Z^d$ implies that $A$ is a Hilbert basis, and then we can
apply again Theorem~\ref{thm:equivalences} (v) implies (ii) to get:

\begin{corollary} \textup{(Sturmfels) \cite[Corollary 8.9]{stu}}
\label{thm:sturm} Let $A \in \Z^{d \times n}$ and let
${\bf c} \in \R^n$ be generic with respect to $A$.
Then $\Delta_{\bf c}$ is a unimodular triangulation
if and only if the toric initial ideal $\textup{in}_{\bf c}(I_A)$
is generated by square-free monomials.
\end{corollary}

Theorem~\ref{thm:equivalences} is the result of generalizing
Sturmfels' above result to arbitrary TDI systems. Still concerning
generic ${\bf c}$ it is worth to note the following result of Conti
and Traverso which provides another connection between integer
linear programming and Gr\"obner bases. Here we think of an element
${\bf x}^{{\bf v}^+} - {\bf x}^{{\bf v}^-}$ as a vector ${\bf v}^+ -
{\bf v}^-$.
\begin{proposition} \textup{(Conti-Traverso) \cite{ct} --
see \cite[Lemma 3]{methods} } \label{lem:ct}
If $\textup{IP}_{A,{\bf c}}({\bf b})$ has a unique optimal solution
for every ${\bf b} \in \N A$ then the reduced Gr{\" o}bner basis
is a minimal test-set for the family of integer programs
$\textup{IP}_{A,{\bf c}}$.
\end{proposition}

This proposition means for us that in the generic case the following
(vi) can be added to Theorem~\ref{thm:equivalences}: (vi) {\em The
initial terms in the reduced Gr{\" o}bner basis are square-free.} In
particular, in the generic case of condition (iii) of
Theorem~\ref{thm:equivalences} the unique inclusionwise minimal
test-set is defined by the reduced Gr\"obner basis, which, by (vi)
has only square-free terms initial terms. Even though
Theorem~\ref{thm:equivalences} concerns general TDI systems and
could be proved in elementary means, it was highly stimulated by
the above results concerning the generic case.

As is typically the case in combinatorial optimization, the cost
vector ${\bf c}$ is not generic for $A$. However, there may be
cases where one can slightly perturb ${\bf c}$ to another cost
vector that is generic for $A$, and where the TDI property for
${\bf y}A \leq {\bf c}$ can be more easily studied when ${\bf c}$
is perturbed. More precisely, from the implication
``(ii) implies (i)''  we immediately get the following:

\begin{proposition} \label{rem:ref}
If ${\bf c}, {\bf c}^\prime \in \R^n$ are such that $\Delta_{{\bf
c}^\prime}$ of $A$ is a refinement of $\Delta_{\bf c}$ of $A$, where
$A_\sigma$ is a Hilbert basis for all $\sigma\in\Delta_{{\bf
c}^\prime}$, and in particular if $\Delta_{{\bf c}^\prime}$ is a
unimodular triangulation of $A$, then ${\bf y}A \leq {\bf c}$ is
TDI.
\end{proposition}

\begin{remark} Having a regular unimodular refinement
$\Delta_{\bf c^\prime}$ of $\Delta_{\bf c}$ amounts to providing
an integer point on the face of each $P_{\bf b}$, for all
${\bf b} \in \textup{cone}(A) \cap \Z^d$, that is minimized by
${\bf c}$. This integer point is the vertex of $P_{\bf b}$
minimized by ${\bf c^\prime}$ and so having a regular unimodular
refinement provides an integer point on each $P_{\bf b}$ in a
uniform manner, dictated by ${\bf c^\prime}$.
\end{remark}

In the rest of the paper one of our favorite themes will be  the use
of Proposition~\ref{rem:ref}. Clearly, the unimodular triangulation
does not even need to be regular: in fact, a unimodular cover of
each of the cells of $\Delta_{\bf c}$ suffices for verifying that
${\bf y}A \leq {\bf c}$ is TDI. However, we are interested in cases
when the converse of Proposition~\ref{rem:ref} is true.  In general
it is not. It is not even true that a Hilbert basis has a unimodular
partition or a unimodular covering \cite{BRUNS} and this counterexample 
inspires two more remarks. First, it cannot be expected
that the equivalence of (i) and (v) can be reduced to Sturmfels'
generic case. Secondly, it should be appreciated that
the converse of this remark does hold in the important set packing
special case, as we will see in the next section.

\section{Set Packing}\label{sec:perfect}

Let a set packing problem be defined with a matrix $A$ and vector
${\bf c}$, and recall  ${\bf c}:=({\bf 1}, {\bf 0}) \in \R^{n}$,
where the last $d$ entries of ${\bf c}$ are $0$. If the set packing
polytope $Q_{\bf c}$ has integer vertices then the matrix $A$ and
the polytope $Q_{\bf c}$ are said to be {\em perfect.} (We will not
use the well-known equivalence of this definition with the integer
values of optima: this will follow.) Lov\'asz' (weak) perfect graph
theorem \cite{lovasz} is equivalent to: {\em the matrix A defining a
set packing polytope is perfect if and only if its first $(n-d)$
columns form the incidence vectors (indexed by the vertices) of the
inclusionwise maximal complete subgraphs of a perfect graph}.

A polyhedral proof of the perfect graph theorem can be split into
two parts: Lov\'asz' {\em replication lemma} \cite{lovasz} and
Fulkerson's  {\em pluperfect graph theorem} \cite{fulkerson}. The
latter states roughly that a set packing polytope with integer
vertices is described by a TDI system of linear inequalities. In
this section we restate Fulkerson's result in a sharper form: there
is a unimodular regular triangulation that refines the normal fan of
any integral set packing polytope. We essentially repeat Fulkerson's
proof, completing it with a part that shows unimodularity along the
lines of the proof of \cite[Theorem 3.1]{sebo}. The following
theorem contains the weak perfect graph theorem and endows it with
an additional geometric feature.

\begin{theorem} \label{thm:perfect}
Let $Q_{\bf c}$ be a set packing polytope defined by $A$ and ${\bf
c}$. Then there exists a vector $\varepsilon \in \R^n$ such that
${\bf c}^\prime :=({\bf 1}, {\bf 0}) + \varepsilon$ defines a
regular triangulation $\Delta_{{\bf c}^\prime}$  refining
$\Delta_{\bf c}$, and this triangulation is unimodular if and only
if $Q_{\bf c}$ is perfect.
\end{theorem}

We do not claim that the following proof of this theorem is novel.
All essential ingredients except unimodularity are already included
in the proof of Fulkerson's pluperfect graph theorem
\cite{fulkerson}.  Fulkerson's proof suggests a greedy way of taking
active rows with an integer coefficient (see below); this is often
exploited to prove that some particular systems are TDI, let us only
cite two papers the closest to ours, Cook, Fonlupt \& Schrijver
\cite{cfs} and Chandrasekaran \& Tamir \cite{ct}. The latter paper
extensively used {\em lexicographically best} solutions, which is an
important tool in linear programming theory. This idea was used in
\cite{sebo} to prove that the normal fan has a unimodular
refinement. This same lexicographic perturbation is accounted for by
the vector $\varepsilon$ of Theorem~\ref{thm:perfect}, showing that
the unimodular refinement is regular. This motivated the following
problem, which contains the perfectness test (detecting perfection in a graph):

\begin{problem}\label{prob:unimod} \cite{paolo}
Given a $d\times n$ integer matrix $A$ and an $n$ dimensional
integer vector ${\bf c}$, decide in polynomial time whether the normal
fan of $Q_{\bf c}$ consists only of unimodular cones. Equivalently, can
it be decided in polynomial time that $Q_{\bf c}$ is non-degenerate,
and the determinant of $A_\sigma$ is $\pm 1$ for all $\sigma \in
\Delta_{\bf c}$.
\end{problem}

Motivated by the perfectness test, the following problem might still be
polynomially solvable.

\begin{problem}\label{pb:pert}
Given a $d\times n$ integer matrix $A$ and an $n$ dimensional
integer vector ${\bf c}$, decide in polynomial time whether
$\Delta_{\bf c}$ has a unimodular refinement $\Delta_{\bf c'}$.
\end{problem}

In the set packing case not only is the perfectness test more
efficient, but also the linear progams and their duals can be solved in
polynomial time (even if the algorithm uses nonlinear optimization
and the ellipsoid method). Could this also be true in general ?

\begin{problem}\label{pb:sol}
Given a $d\times n$ integer matrix $A$ and an $n$ dimensional
integer vector ${\bf c}$ as input, is there a polynomial algorithm
that finds at least one of the following as output: a solution to
$\textup{IP}_{A,{\bf c}}({\bf b})$, or a ``NO'' answer to
Problem~\ref{pb:pert} (or to Problem~\ref{prob:unimod}).
\end{problem}

All these problems can be solved in polynomial time in the set
packing case: Problem~\ref{pb:pert} according to
Theorem~\ref{thm:perfect} for the moment only by using
\cite{recognition}; Problem~\ref{pb:sol} for the moment only by
\cite{gls}. It remains an interesting question whether there are
more efficient and conceptually simpler solutions to these problems.

We now prepare the proof of Theorem~\ref{thm:perfect}. It is a last
step in a sharpening series of observations all having essentially
the same proof.
We begin with the proof of Fulkerson's pluperfect graph
theorem which will indicate what the ${\bf c}'$ of
Theorem~\ref{thm:perfect} should be, and then finish by showing that
$\Delta_{{\bf c}'}$ is a unimodular triangulation.

Assume that $A$ is a perfect matrix for the remainder of this
section and that ${\bf c} = ({\bf 1},{\bf 0})$ as before.
For all ${\bf b} \in \Z^d$ and $ i \in [n]$ let
$$
\lambda_{{\bf c},i}({\bf b}) := \textup{max} \{ x_i : \, {\bf x} \,
\textup{ is an optimal solution of $\textup{LP}_{A,{\bf c}}({\bf
b})$ }\}.
$$
That is, $\lambda_{{\bf c},i}({\bf b})$ is the largest value of
$x_i$ such that ${\bf c} \cdot {\bf x}$ is minimum under ${\bf x}
\in P_{\bf b}$.

\begin{remark}
If $\sigma$ is the minimal cell of
$\Delta_{\bf c}$ such a ${\bf b} \in \textup{cone}(A_\sigma)$, then
$ {\bf b} - \lambda_{{\bf c},i}({\bf b}) {\bf a}_i \in
\textup{cone}(A_{\sigma^\prime}) $ where $\sigma^\prime \in
\Delta_{\bf c}$, $\sigma^\prime \subseteq \sigma$ and the dimension
of $\textup{cone}(A_{\sigma^\prime})$ is strictly smaller than that
of $\textup{cone}(A_\sigma)$. Furthermore, ${\bf b} - \lambda {\bf
a}_i \notin \textup{cone}(A_\sigma)$ if $\lambda > \lambda_{{\bf
c},i}({\bf b})$.
\end{remark}

For all ${\bf b} \in \Z^d$ we show that $\lambda_{{\bf c},i}({\bf b})$ is an
integer for every $i = 1,\ldots,n$. This is the heart of Fulkerson's
pluperfect graph theorem \cite[Theorem 4.1]{fulkerson}. We state this in the
following lemma in a way that is most useful for our needs.
Denote the common optimal value of $\textup{LP}_{A,{\bf c}}({\bf
b})$ and $\textup{DP}_{A,{\bf c}}({\bf b})$ by $\gamma_{\bf c}({\bf
b})$. Note that $\gamma_{\bf c}$ is a monotone increasing function
in all of the coordinates.

\begin{lemma}  \label{lem:allstrong}
Suppose $\gamma_{{\bf c}}({\bf b})\in\Z$ for all  ${\bf b} \in \Z^d$.
If ${\bf x}$ is an optimal solution to
$\textup{LP}_{A,{\bf c}}({\bf b})$ with $x_l \neq 0$ for some
$1 \leq l \leq n$, then there exists ${\bf x}^*$ also optimal
for $\textup{LP}_{A,{\bf c}}({\bf b})$, such that $x_l^* \geq 1$.
\end{lemma}

Note that this lemma implies the integrality of $\lambda :=
\lambda_{{\bf c},l}({\bf b})$ for all $l=1,\ldots, n$: if $\lambda$
were not an integer then setting ${\bf b}^\prime := {\bf b} -
\lfloor \lambda \rfloor {\bf a}_l$ we have $\lambda_{{\bf c},l}({\bf
b}^\prime) = \{ \lambda \}$ where $ 0 \le \{\lambda\} := \lambda -
\lfloor \lambda \rfloor <1$, contradicting
Lemma~\ref{lem:allstrong}.

\begin{proof}
Suppose ${\bf x} \in P_{\bf b}$ with
${\bf c} \cdot {\bf x} = \gamma({\bf b})$ and $x_l > 0$ for some
$1 \leq l \leq n$. We can assume that $x_l < 1$ since otherwise
${\bf x}^* := {\bf x}$. We have two cases: either
$1 \leq l \leq n-d$ or $n-d+1 \le l \le n$.

If $n-d+1 \le l \le n$ then ${\bf a}_l= -{\bf e}_{l-(n-d)} \in \R^d$
and $c_l = 0$. In this case, we have $\gamma_{\bf c}({\bf b}) =
\gamma_{\bf c}({\bf b}+ x_l {\bf e}_{l-(n-d)})$ because replacing
$x_l$ by $0$ in ${\bf x}$ we get a solution of the same objective
value for the right hand side ${\bf b}+ x_l {\bf e}_{l-(n-d)}$ which
gives $\gamma_{\bf c}({\bf b}) \geq \gamma_{\bf c}({\bf b}+ x_l {\bf
e}_{l-(n-d)})$. The reverse inequality follows from the
(coordinate-wise) monotonicity of $\gamma_{\bf c}$. But then
$$
\gamma_{\bf c}({\bf b} + {\bf e}_{l-(n-d)})
\le \gamma_{\bf c}({\bf b} + x_l {\bf e}_{l-(n-d)}) +1 - x_l
= \gamma_{\bf c}({\bf b}) + 1 - x_l,
$$
and since $\gamma_{\bf c}({\bf b}+ {\bf e}_{l-(n-d)})$ is integer
and  $0 < 1-x_l <1$, we conclude that $\gamma_{\bf c}({\bf b}+ {\bf
e}_{l-(n-d)}) = \gamma_{\bf c}({\bf b})$.

So for any optimal
${\bf x}^\prime \in P_{{\bf b}+ {\bf e}_{l-(n-d)}}$ where
${\bf c} \cdot {\bf x}^\prime = \gamma_{\bf c}({\bf b})$,
letting ${\bf x}^*:= {\bf x}^\prime + {\bf e}_{l-(n-d)} \in P_{\bf b}$
we have ${\bf c} \cdot {\bf x}^* \le \gamma_{\bf c}({\bf b})$
and so ${\bf x}^*$ is optimal and $x^*_l \ge 1$.

Suppose now $1 \leq l \leq n-d$. Replacing $x_l$ in ${\bf x}$ by $0$
we get a point in $P_{{\bf b}-x_l {\bf a}_l}$. This point has
objective value ${\bf c} \cdot {\bf x} - x_l < {\bf c} \cdot {\bf x}
= \gamma_{\bf c}({\bf b})$, and so we have by monotonicity
$$
\gamma ({\bf b}-{\bf a}_l) \le \gamma ({\bf b}-x_l{\bf a}_l) <
\gamma ({\bf b}).
$$
Since the left and right hand sides are both integer values then
$\gamma ({\bf b}-{\bf a}_l) \le \gamma ({\bf b}) - 1$.
Letting ${\bf x}^* := {\bf x}^\prime + {\bf e}_l \in P_{\bf b}$
we have
$
{\bf c} \cdot {\bf x}^*
\le \gamma_{\bf c}({\bf b})-1 +1
= \gamma_{\bf c}({\bf b}) $,
so $x^*$ is optimal, and  $x^*_l\ge 1$.
\end{proof}

Let us now define the appropriate ${\bf c}'$ for the theorem,
depending only on ${\bf c}$. Define ${\bf c}^\prime := {\bf c} +
\varepsilon \in \R^n$ where $\varepsilon_i:= - (1/d^{d+2})^i$ for
each $i=1,\ldots , n$. Note that the absolute value of the
determinant of a $d \times d$  $\{-1, 0, 1 \}$-matrix cannot exceed $d^d$.
It follows, by Cramer's rule, that the coefficients of linear
dependencies between the columns of $A$ are at most $d^d$ in
absolute value, and then the sum of absolute values of the
coefficients between two solutions of an equation
$A{\bf x} = {\bf b}$ for any $ {\bf b} \in \R^d$ can differ by at
most a factor of $d^{d+2}$. After this observation, the following lemma
is straightforward to verify

\begin{lemma} \label{lem:easy}
\begin{itemize}
\item[(i)] Any optimal solution to
$\textup{LP}_{A, {\bf c'}}({\bf b})$ is also optimal for
$\textup{LP}_{A, {\bf c}}({\bf b})$.
\item[(ii)] If ${\bf x}'$ and ${\bf x}''$ are both
optimal solutions to $\textup{LP}_{A, {\bf c}}({\bf b})$ then ${\bf
x}'$ is {\em lexicographically bigger} than ${\bf x}''$ (that is,
the first non-zero coordinate of ${\bf x}' - {\bf x}''$ is positive)
if and only if ${\bf c}' \cdot {\bf x}' < {\bf c}' \cdot {\bf x}''$.
\end{itemize}
\end{lemma}

Statement (i) of the lemma means that $\Delta_{\bf c'}$ refines $\Delta_{\bf c}$,
and statement (ii) means that an optimal solution to
$\textup{LP}_{A, {\bf c'}}({\bf b})$
is constructed by defining ${\bf b}^0 := {\bf b}$ and recursively
$$
x_i:=\lambda_{{\bf c},i}({\bf b}^{i-1}), \,\,\, {\bf b}^i := {\bf
b}^{i-1} - x_i{\bf a}_i \, \, \, \textup{for all} \, i=1,\ldots, n.
$$
Furthermore, this optimum is unique and it follows that
$\Delta_{\bf c'}$ is a triangulation. We are now ready to
prove Theorem~\ref{thm:perfect}.

\vspace{.2cm}

\noindent {\em Proof of Theorem~\ref{thm:perfect}}. The necessity of
the condition is straightforward: each vertex ${\bf y} \in Q_{\bf c}$
satisfies the linear system (consisting of $d$ equations in $d$ unknowns)
of the form ${\bf y} A_{\sigma'} = {\bf c}_{\sigma'}$ where $\sigma'$ is a cell
of $\Delta_{{\bf c}'}$, ${\bf c}_{\sigma'}$ is the subvector of ${\bf c}$
indexed by $\sigma'$.  Since the determinant of $A_{\sigma'}$ is $\pm 1$,
${\bf y}$ must be an integer vector because of Cramer's rule.

Conversely, we will prove the assertion supposing only that
$\gamma_{\bf c}({\bf b})$ is integer for all ${\bf b} \in \Z^d$
and applying Lemma~\ref{lem:allstrong}. Note that, by the already
proven easy direction, we will have proved from this weaker
statement that $Q_{\bf c}$ is perfect.
Without loss of generality, suppose that ${\bf b} \in \Z^d$ cannot be
generated by less than $d$ columns of $A$. That is, the minimal cell
$\sigma$ of $\Delta_{\bf c}$ such that ${\bf b} \in \textup{cone}(A_\sigma)$
is a maximal cell of $\Delta_{\bf c}$. That is,
$\textup{cone}(A_\sigma)$ is $d$-dimensional.

Because of Lemma~\ref{lem:easy}(i), an optimal solution to
$\textup{LP}_{A,{\bf c'}}({\bf b})$ will have support in $\sigma$ and
Lemma~\ref{lem:easy}(ii) implies that such an optimal
solution is constructed as follows:
Let $s_1 := \textup{min} \{ i \, : \, i \in \sigma \}$ and $x_{s_1}
:= \lambda_{{\bf c},s_1}({\bf b})$. Recursively, for $j = 2, \ldots,
d$ let $s_j$ be the smallest element in $\sigma$ indexing a column
of $A$ on the minimal face of $\textup{cone}(A_\sigma)$ containing
$$
{\bf b} -\sum_{i=1}^{j-1}x_{s_i}{\bf a}_{s_i}.
$$

Since ${\bf b}$ is in the interior of $\textup{cone}(A_\sigma)$ then
$x_{s_i}>0$ for each $i=1, ...,d$, and by Lemma~\ref{lem:allstrong},
these $d$ $x_{s_i}$'s are integer. Moreover, since the dimension of
$\textup{cone}(A_{\sigma \backslash \{s_1, \ldots, s_i \}})$ is
strictly decreasing as  $i=2, \ldots, d$  progresses, then
$
{\bf b} \, - \, \sum_{i=1}^{d}x_{s_i}{\bf a}_{s_i} \, = \, 0
$,
and setting
$U := \{ s_1,\ldots, s_d \} \subseteq \sigma$, we have that the
columns of $A_U$ are linearly independent. Note that $U$ is a cell
of $\Delta_{\bf c'}$ and by Lemma~\ref{lem:easy}(ii) every maximal
cell of $\Delta_{\bf c'}$ arises in this fashion. We show that the
matrix $A_U$ has determinant $\pm 1$.

Suppose not. Then the inverse of the matrix $A_U$ is non-integer,
and from the matrix equation $(A_U)^{-1}A_U \, = \, I_d$ we
see that there exists a unit vector ${\bf e}_j \in \R^d$ which is a
noninteger combination of columns in $A_U$: $ \sum_{i=1}^d x_{s_i}
{\bf a}_{s_i} = {\bf e}_j $. Let ${\bf z}$ be the vector
$$
{\bf z} : = \sum_{i=1}^d \{x_{s_i}\} {\bf a}_{s_i}.
$$
Clearly, ${\bf z} \in \textup{cone} (A_U)$ and furthermore
${\bf z} \in \Z^d$ since it differs from ${\bf e}_j$ by an integer
combination of the columns of $A_U$. So Lemma~\ref{lem:allstrong}
can be applied to ${\bf b} := {\bf z}$: letting $l := \textup{min}
\{ i  \, : \, \{ x_{s_i}\} \neq 0 \}$ we see that $\lambda_{{\bf
c},s_l}({\bf z})=\{x_{s_i}\} < 1$ contradicting
Lemma~\ref{lem:allstrong}. Hence both $A_U$ and $(A_U)^{-1}$ are
integer,  their determinant is $\pm 1$; since $A_U$ was an arbitrary
maximal cell of $\Delta_{\bf c'}$, we conclude that $\Delta_{\bf
c'}$ is unimodular. \qed

\vspace{.2cm}

The argument concerning the inverse matrix replaces the use of
parallelepipeds (compare with \cite[proof of Theorem 3.1]{sebo}).
Note that all the numbers in the definition of ${\bf c'}$ are at most
$d^{d^2}$, so they have a polynomial number of digits: the perturbed
problem has polynomial size in terms of the original one, reducing
the perfectness test to Problem~\ref{prob:unimod}.

\begin{remark}
Note that the $\varepsilon$ of the perturbed vector
${\bf c}^\prime$ in Theorem~\ref{thm:perfect} made
no prescribed order on the columns of $A$. These regular
refinements are known as {\em pulling}
refinements ~\cite{lee} and so Theorem~\ref{thm:perfect}
complements the result of Ohsugi and Hibi \cite{oh}
regarding pulling refinements, of not the normal fan
of $Q_{\bf c}$ but, of the polytope $Q_{\bf c}$ itself.
See \cite{lee} for more background on triangulations of polytopes.
\end{remark}

\section{Computation}\label{sec:computation}

In this section we wish to provide an illustration of how the
results presented in this work lead to practical algorithms. We
argue that the computational algebra methods for detecting the TDI
property
 can be especially efficient in practice when there is a
generic perturbation of the system in the sense of
Proposition~\ref{rem:ref}. We wonder if this practical efficiency is
in some way related to the detection of perfection being in P
\cite{recognition} ?  Could  the perfection recognition problem be
solved with a polynomial algorithm based on such geometric ideas ?
By analogy to the perfectness test, could a unimodular perturbation
be found in polynomial time (Problem~\ref{pb:pert}) ?

To show the ideas, we focus on one example of a set packing problem
for a chordal graph $G$ with $6$ maximal cliques on $10$ vertices. The coefficient
matrix $A$ has $10$ rows and $16$ columns and the cost vector ${\bf
c}$ has $16$ entries; the first $6$ being equal to $1$, the last
$10$ of which equal $0$ (see the Figure).

\begin{figure}
\centering
\input{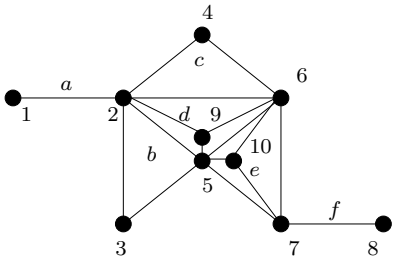}
\caption{A chordal graph G with 6 maximal cliques on 10 vertices.}
   \label{fig:padberg}
\end{figure}

There are essentially five different ways in which we
can detect the TDI property in the system ${\bf y}A \leq {\bf c}$
coming from the graph $G$, in finite time, and in polynomial time
in fix dimension.
The first two use Theorem~\ref{thm:perfect}:
Indeed, one can test (even if the time will depend
exponentially on the dimension) if the triangulation $\Delta_{\bf
c^\prime}$ is unimodular. The second possibility is to test if the
(monomial) toric initial ideal $\textup{in}_{\bf c^\prime}(I_A)$ is
generated by squarefree monomials, by computing a Gr\"obner basis.

There are two other ways using Theorem~\ref{thm:equivalences}.
Ignoring the generic perturbation arising from
Theorem~\ref{thm:perfect}, we can simply study the original set
packing system ${\bf y}A \leq {\bf c}$: a third algorithm can be
based on statement (iii) of Theorem~\ref{thm:equivalences}, or more
precisely, on Corollary~\ref{thm:cls}. A fourth would be to use
statement (v) of Theorem~\ref{thm:equivalences}.

A fifth possibility for the set packing special case is to use the
polynomial algorithm of Chudnovsky, Cornu{\' e}jols, Liu, Seymour
and Vu\v{s}kovi\'c \cite{recognition} which (and actually its early
predecessors for testing whether a graph is chordal) are certainly
the best on this example and any chordal graph, since it provides an
immediate positive answer through clique separation. However, the 
example will only contrast the methods expounded upon here.

\noindent {\bf The
triangulation $\Delta_{\bf c^\prime}$:} Using {\tt polymake}
\cite{polymake} we can compute the vertices of the simple polytope
$Q_{\bf c^\prime}$. We find that there are $288$ vertices in all and
for each of these vertices, we can find their $10$ active facets.
Finally, for each vertex the coordinates of the active facets give
rise to a $10 \times 10$ matrix and we would need to check that each
of these $288$ square matrices have determinant $\pm 1$.

\noindent
{\bf The initial ideal $\textup{in}_{\bf c^\prime}(I_A)$:}
The toric ideal $I_A$ lives in ${\bf k}[a, \ldots, f, v_1, \ldots, v_{10}]$
where $a,\ldots,f$
correspond to the maximal cliques of $G$ (the first 6 columns of $A$)
and where $v_1,\ldots,v_{10}$ correspond to the vertices of $G$
(the ordered columns of $-I_{10}$, the last $10$ columns of $A$) as
before. The toric ideal for set packing matrices $A$ has a very
simple generating set: there is one generator for every maximal
clique in the graph, all of the form
$z v_{i_1}v_{i_2} \cdots v_{i_r} - 1$, where $z$ is the variable
corresponding to the clique and $\{i_1, i_2, \ldots, i_r \}$ are
the vertices of that clique. See \cite[\S 2.2]{oshea} for justification.

This simple generating set enables
quick computation of the toric initial ideals of $A$. Using
{\tt Macaulay 2} \cite{M2}, we computed
$ \textup{in}_{\bf c^\prime}(I_A) \, = \,
\langle \,
fv_7v_8, \,  ev_5v_6v_{10}, \, dv_2v_9, \,  cv_4, \,  bv_3, \,  av_1 \,
\rangle
$
in less than one second on a standard desktop.
The monomial toric initial ideal $\textup{in}_{\bf c^\prime}(I_A)$
has 6 generators and all are squarefree. Not only was this initial
ideal computed quickly but furthermore, because of the small number
of generators, we see that the TDI property is presented in a
highly compact manner, unlike the triangulation $\Delta_{{\bf c}^\prime}$
above.

\noindent
{\bf The construction of ${\mathcal T}_{A,{\bf c}}$:}
To construct ${\mathcal T}_{A,{\bf c}}$ we would first
need to find all sets of size $10$ that are not cells
in $\Delta_{\bf c}$. There are a total of
${16}\choose{10}$ $= 8008$ candidate
sets of size $10$. We would then have to eliminate all
$10$-sets that are contained in some cell of $\Delta_{\bf c}$.
Again, using {\tt polymake} \cite{polymake}, we can find
all $101$ maximal cells of $\Delta_{{\bf c}}$.
One such cell is the $14$-set
$\sigma := [16] \backslash \{8,13\}$ and
so all $10$-sets that are contained
in $\sigma$ must be removed. Similarly, using the other $100$
cells of $\Delta_{\bf c}$ we could find all $10$-sets that are
non-cells.
We would then need to implement the remaining steps outlined
in the algorithm of Corollary~\ref{thm:cls}, a computationally
challenging process that would need to be carried out for each
of the non-cells. This task does not compare favourably to the
computation of $\textup{in}_{\bf c^\prime}(I_A)$ above.

\noindent
{\bf The monomial ideal
$\textup{mono}(\textup{in}_{\bf c^\prime}(I_A))$:}
There is an algorithm \cite[Algorithm 4.4.2]{sst} for computing
the generators of the monomial ideal
$\textup{mono}(\textup{in}_{\bf c}(I_A))$ which involves four steps,
all implemented in {\tt Macaulay 2} \cite{M2}.
The first is finding a generating set for the toric initial ideal
$\textup{in}_{\bf c}(I_A)$; one way to do this is to find
a Gr{\" o}bner basis of $I_A$ with respect to ${\bf c}$. For our
example, one generating set is
$\{ f v_7v_8, \, e v_5 v_6 v_{10} - f v_8, \, d v_2 v_9 - e v_7 v_{10}, \,
c v_4 - d v_5 v_9, \,  b v_3 - d v_6 v_9, \, a v_1 - d v_5 v_6 v_9 \}$
which consists of $1$ monomial and
$5$ binomials in the polynomial ring
${\bf k}[a, \ldots, f, v_1, \ldots, v_{10}]$.
Next, we must carry out a {\em multi-homogenization} procedure
on this generating set to form the multi-homogenized ideal
$\textup{in}_{\bf c}(I_A)^\textup{homo}$ in the
polynomial ring
${\bf k}[a, \ldots, f, v_1, \ldots, v_{10},
A,\ldots, F, V_1, \ldots, V_{10} ]$
consisting of $32$ variables where the $16$ new homogenizing
variables are $A,\ldots, F, V_1, \ldots, V_{10}$. Finding
generating sets for $\textup{in}_{\bf c}(I_A)$ and
$\textup{in}_{\bf c}(I_A)^\textup{homo}$ can be quickly done
using {\tt Macaulay~2}.

The third step is to find a Gr{\" o}bner basis of
$\textup{in}_{\bf c}(I_A)^\textup{homo}$ with respect to
an elimination order.
The fourth and final step would be to extract the monomials in
the Gr{\" o}bner basis of step 3 -- this set of monomials generate the
monomial ideal:
$
\textup{mono}(\textup{in}_{\bf c}(I_A)) \, = \,
\langle \,
av_1v_2, \, bcv_2v_3v_4, \, bv_2v_3v_5, \, cv_2v_4v_6, \, bfv_2v_3v_8, \, cfv_2v_4v_8, \, fv_7v_8, \, dv_2v_5v_6v_9, \, dfv_2v_8v_9, \,
aev_1v_7v_{10},
$ \\
$
bcev_3v_4v_7v_{10}, \, bev_3v_5v_7v_{10}, \, cev_4v_6v_7v_{10}, \, ev_5v_6v_7v_{10}
\, \rangle.
$
Like the computation of $\textup{in}_{{\bf c}^\prime}(I_A)$ above, this set of
$14$ generators was done in less than a second on a standard desktop and
from it we can see the TDI propertry presented in a compact manner. Note
that, from the equivalence of (iii) and (v) in Theorem~\ref{thm:equivalences},
we also get the $14$ $\kappa$'s and $\tau$'s of ${\mathcal T}_{A,{\bf c}}$ from
these generators.

We hope that this pedagogical example exhibits some of the utility
of computational algebra methods in integer programming. In larger
examples, like Padberg's {\em windmill} \cite[last figure]{padberg}
(a graph with $21$ maximal cliques and $20$ vertices), the
computation of $\textup{mono}(\textup{in}_{\bf c}(I_A))$ cannot be
carried out on a standard desktop but computing $\textup{in}_{\bf
c^\prime}(I_A)$ can be carried out in less than a second. However,
it should be noted, that for both our example above and for
Padberg's windmill, it is easy to observe that the graph is perfect
via clique separation, a classical simple operation that preserves
perfectness.

The ability to only attain $\textup{in}_{\bf c^\prime}(I_A)$ in the
Padberg example was typical of other examples of large set packing
problems. Could this or some other ``geometric'' (polyhedral) method
ever handle perfectness -- a notion that can be defined in purely
polyhedral terms -- efficiently (in polynomial time) and in a less
technical way  ?  In view of the NP-completeness of the general
problem \cite{chinese}, \cite{juli}, unimodular covering may provide
a distinguishing clue for this $0-1$ special case.

\section*{Acknowledgments}
The authors wish to thank Rekha Thomas for her valuable input and
suggestions. Some more work related to the results of this article,
including the computational experimentation, can be found in the
first author's Ph.D. dissertation at the University of Washington.
Thanks are also due to the developers of the computational packages
used in this work.

\end{document}